\documentclass[11pt]{article}
\usepackage{amsmath}
\usepackage[latin1]{inputenc}
\usepackage{amsthm}
\usepackage{amssymb}
\usepackage{amscd}
\usepackage{enumerate}
\usepackage{mathtools}
\newtheorem{teo}{\bf Theorem}

\newtheorem{cor}[teo]{Corollary}
\newtheorem{lema}[teo]{Lemma}
\theoremstyle{remark}
\newtheorem{remark}{\bf Remark}
\DeclareMathOperator*{\esup}{ess\,sup}
\DeclareMathOperator*{\einf}{ess\,inf}
\renewcommand{\thefootnote}{}
\begin{document}

\date{}

\title{\bf Discrete convolution operators and  Riesz systems generated by actions of abelian groups}

\author{
Gerardo P\'erez-Villal\'on
}

\maketitle

{\def\thefootnote{}
}

{\def\thefootnote{}
\footnotetext{{\it E-mail:} gerardo.perez@upm.es}}
\footnotetext{Departamento de Matem\'atica Aplicada a las TIC, UPM,
 Nikola Tesla, s/n
28031 Madrid, Spain. }
{\def\thefootnote{}

\begin{abstract}
We study the bounded endomorphisms of $\ell_{N}^2(G)=\ell^2(G)\times \dots  \times\ell^2(G)$ that commute with translations, where $G$ is a discrete abelian group. It is shown that they form a C*-algebra isomorphic to the C*-algebra of $N\times N$ matrices with entries in $L^\infty(\widehat{G})$, where $\widehat{G}$ is the dual space of $G$.
Characterizations of  when these endomorphisms are  invertible, and  expressions for their norms and for the norms of their inverses, are given. These results allow us to
study Riesz systems that arise from the action of $ G $ on a finite set of elements of a Hilbert space.
\end{abstract}

\section{Introduction}

  Let $G$ be a discrete abelian group.
The first aim of this work is to study operators of the type
\[
\mathcal{A}:\ell_{N}^2(G)\rightarrow \ell_{M}^2(G),\qquad \mathcal{A}(\mathbf{x})=A\ast \mathbf{x}=\sum_{g\in G} A(g) \, \mathbf{x}(\cdot-g),
\]
where $A\in \mathcal{M}_{_{M\times N}}(\ell^2(G))$ is a $M\times N$ matrix whose
entries are elements of $\ell^2(G)$ and $A\ast \mathbf{x}$ is the convolution of the matrix $A$ with the vector
$\mathbf{x}\in \ell^2_{_{N}}(G)=\ell^2(G)\times\dots\times \ell^2(G)$  ($N$ times).

 The bounded convolution operators of this type can also be described as  those bounded linear operators that commute with translations (see Theorem \ref{bessel2}). When $ G = \mathbb {Z} ^ d $ they are called in discrete signal processing,   where they are widely used, Linear Time Invariant (LTI) Multi-Input Multi-Output (MIMO) transformations. See for example \cite{kailath:80}.

In many situations, operators between some spaces of functions or measures, on  a locally compact abelian group $G$, that commute with translations coincide with
those that can be expressed as a multiplication in the Fourier domain (see Theorem \ref{bessel2} for the discrete vectorial case here considered). They are called multipliers and have been widely studied for escalar functions and measures. See for example \cite{larsen:71}, in particular the closest result to this work   \cite[Theorem 4.3.1]{larsen:71},
 where characterizations of the multipliers in $\mathcal{B}(L^2(G))$ are given for a general locally compact abelian group $G$.

 We give, in Seccion 3, characterizations of the multipliers from $\ell_{N}^2(G)$ into
$\ell_{M}^2(G)$, where $G$ is a discrete abelian group. We study some of the characteristics of these discrete vectorial convolution operators. Special attention is devoted to the convolution operators in $\mathcal{B}(\ell^2_N(G))$, the space of bounded endomorphism of $\ell^2_N(G)$. It is proved that they form a C*-subalgebra of $\mathcal{B}(\ell^2_N(G))$ which is isomorphic to $\mathcal{M}_N\big(L^\infty(\widehat{G})\big)$, the C*-algebra of the $N\times N$ matrices with entries in $L^\infty(\widehat{G})$,
 where $\widehat{G}$ is the dual space of $G$. For instance, the set of convolution operators of $\mathcal{B}(\ell^2_N(\mathbb{Z}))$ form a C*-algebra  isomorphic to $\mathcal{M}_N\big(L^\infty(\mathbb{T})\big)$ where $\mathbb{T}=\{z\in \mathbb{C}: |z|=1\}$, or those in $\mathcal{B}(\ell^2_N(\mathbb{Z}^2))$ form a C*-algebra  isomorphic to $\mathcal{M}_N\big(L^\infty(\mathbb{T}^2)\big)$.

 By means of this C*-isomorphism, we characterize the invertible convolution operators, and we give a suitable expression for the norm of the inverse.

%\medskip

These results about discrete convolution operator in this general setting could be useful in future applications, specially in discrete signal processing where not only the discrete group $\mathbb{Z}^d$, but also finite groups  such as  $\mathbb{Z}^d_s$ or direct products as $\mathbb{Z}_s \times \mathbb{Z}_r$ or  $\mathbb{Z}_s \times \mathbb{Z}^d$, are often used. The second part of this article, shows that by means of these convolution operators, an interesting generalization of some relevant results about shift-invariant systems can be obtained. In  reference \cite{garcia:19b}, it is showed that they are also useful in order to obtain a regular sampling theory in a very general context.

\medskip
 The second aim of this work is to study  Bessel and Riesz systems generated by actions of abelian groups.
The development of wavelet and approximation  theories in different directions has led to the consideration and analysis of various generalizations of the classical shift-invariant spaces in $\mathbb{R}$,
\[
V_\Phi=\Big\{ \,  \sum_{n=1}^N\sum_{\alpha\in \mathbb{Z}} x_{n}(\alpha)\, \varphi_{n}(\cdot-\alpha) \,:\,\, (x_{1},\ldots,x_{N})\in \ell_{N}^2(\mathbb{Z})\, \Big\} \subset L^2(\mathbb{R}),
\]
where $\Phi=\{\varphi_1,\ldots,\varphi_N\}$ denotes the set of generators.
Here, we consider the spaces
\[
V_\Phi=\Big\{ \sum_{n=1}^N \sum_{g\in G} x_{n}(g)\,\, \pi_{g}\varphi_{n} \, :\,\,(x_{1},\ldots,x_{N})\in \ell_{N}^2(G)\,\Big\}\subset \mathcal{H},
\]
where $\mathcal{H}$ is a Hilbert space, $G$ a discrete abelian group, and $\pi$ a unitary representation of $G$ on $\mathcal{H}$.  This generalization of shift-invariant spaces is related with many of the generalizations previously considered, see \cite{aldroubi:96, barbieri:15, barbieri:18, cabrelli:10, hernandez:10} and Section 4 for more details.

Characterizations of when the system $\big\{ \pi_{g}\varphi_{n} \big\}_{n=1,\ldots,N, g\in G}$ is a Riesz basis for $V_\Phi$, of when it is a Bessel sequence of $\mathcal{H}$ and suitable expressions, in the Fourier domain, for the optimal Riesz bounds are provided in Section 4 using the results on convolution operators of the Section 3.

\medskip

\section{Notation and preliminaries}

Throughout the paper we assume that $G$ is a discrete abelian group  (with additive notation), that $\mathcal{H}$ is a separable complex Hilbert space, and we use the following notation:

\begin{enumerate}
\item $\ell^2(G)=\big\{x:G\to \mathbb{C}\, : \|x\|^2_{\ell^2(G)}=\sum_{g\in G}  |x(g)|^2<\infty\big\}$ and
$\ell^2_{_N}(G)=\ell^2(G)\times\ldots\times \ell^2(G)$  ($N$ times).
\item $[\mathbf{x}]_{m}=x_{m}$ \, denotes the $m$-th entry of $\mathbf{x}\in \ell^2_{_{N}}(G)$.
\item $T_{g} \mathbf{x} = \mathbf{x}(\cdot-g)$ denotes the translation operator in $\ell_{_N}^2(G)$ and also in $\ell_{_M}^2(G)$.
\item $\mathcal{M}_{_{M\times N}}(\ell^2(G))$  is the set
of $M\times N$ matrices with entries in $\ell^2(G)$ and  $\mathcal{M}_{_{ N}}(\ell^2(G))=\mathcal{M}_{_{N\times N}}(\ell^2(G))$.
\item $\mathcal{B}\big(\ell^2_{_{N}}(G)\big)$ is the algebra of the bounded endomorphisms of $\ell^2_{_{N}}(G)$.
\item The symbol $\ast$ denotes convolution, namely; $x\ast y=\sum_{g\in G}x(g)y(\cdot-g)$, for  $x,y\in \ell^2(G)$;  $A\ast B=\sum_{g\in G} A(g)B(\cdot-g)$, for $A,B\in \mathcal{M}_{_{ N}}(\ell^2(G))$; $A\ast \mathbf{x}= \sum_{g\in G} A(g)\mathbf{x}(\cdot-g)$, for $A=[a_{m,n}] \in \mathcal{M}_{_{ M\times N}}(\ell^2(G))$ and $\mathbf{x} \in \ell^2_N(G)$,   or equivalently
   $A\ast \mathbf{x}$ is the vector whose $m$-th entry is $[A\ast \mathbf{x}]_{m}=\sum_{n=1}^N a_{m,n}\ast x_{n}$.
\item $\|D\|_{2}$ is the spectral norm of a matrix $D$. The symbols $\lambda_{\min}[C]$ and  $\lambda_{\max}[C]$ denote the minimun and the maximun eigenvalue of a positive semidefinite  matrix $C$.
\item $\pi$ is a unitary representation of the group $G$ on $\mathcal{H}$ i.e. a homomorphism $\pi:G\rightarrow U(\mathcal{H})$, where $U(\mathcal{H})$ is the group of unitary operators of $\mathcal{H}$, that satisfies
$\pi_{g+g'}=\pi_{g}\pi_{g'}$ and $\pi_{-g}=\pi_{g}^{-1}$ for all $g,g'\in G$.
 \item Since $G$ is discrete, its dual space $\widehat{G}$ is compact.  We normalize the Haar measure of $\widehat{G}$ so that $\mu(\widehat{G})=1$ and we define $L^\infty(\widehat{G})$ and $L^2(\widehat{G})$ as usual. Let $\mathcal{F}(x)=\widehat{x}$  denote the Fourier transform of $x$ which is defined by  $\widehat{x}(\xi)=\sum_{g\in G}x(g) \langle -g,\xi \rangle$ for $x$ in $\ell^1(G)$ and it is extended by density to a bijective isometry between  $\ell^2(G)$ and $L^2(\widehat{G})$. See
 \cite[4.5 and 4.6]{folland:95} or  \cite[1.2.7 and 2.2.2]{rudin4} for the most common cases, $G=\mathbb{Z}$, $\mathbb{Z}^d$, $\mathbb{Z}_{s}=\mathbb{Z}/s\mathbb{Z}$, $\mathbb{Z}^d_{s}$, $\mathbb{Z}_{r}\times\mathbb{Z}_{s}$, $\mathbb{Z}^d\times \mathbb{Z}_{r}$,  $\mathbb{Z}\times \mathbb{Z}_{r}\times \mathbb{Z}_{s}$. For instance:
 \[
 \begin{split}
 & \widehat{\mathbb{Z}}\cong\mathbb{T}\quad \text{and}\quad \widehat{x}(z)=\sum_{n\in \mathbb{Z}} x(n)z^{-n},\,\, z\in \mathbb{T}=\{z\in \mathbb{C}: |z|=1\},
 \\& \widehat{\mathbb{Z}^d}\cong \mathbb{T}^d\quad \text{and}\quad \widehat{x}(\mathbf{z})= \sum_{n\in \mathbb{Z}^d} x(n)\mathbf{z}^{-n}, \,\, \mathbf{z}=(z_1,\ldots,z_d)\in \mathbb{T}^d, \, \text{where}\,\, \mathbf{z}^n= z_1^{n_1}\ldots z_d^{n_d},
 \\& \widehat{\mathbb{Z}_s}\cong \mathbb{Z}_s\quad \text{and}\quad \widehat{x}(m)=\sum_{n\in \mathbb{Z}_s } x(n) W_s^{-nm}, \,\, m\in \mathbb{Z}_s, \, \text{where}\,\, W_s=e^{2\pi i/s},
 \\& \widehat{\mathbb{Z}^d\times \mathbb{Z}_s}=\mathbb{T}^d\times\mathbb{Z}_s\quad \text{and}\quad \widehat{x}(\mathbf{z},m)=\sum_{\mathbf{z}\in \mathbb{T}^d, n\in \mathbb{Z}_s } x(\mathbf{z},n)\mathbf{z}^{-n} W_s^{-nm}, \,\, \mathbf{z}\in \mathbb{T}^d,  m\in \mathbb{Z}_s.
 \end{split}
 \]
\item $\mathcal{M}_{_{M\times N}}(L^\infty(\widehat{G}))$  is the set of $M\times N$ matrices with entries in $L^\infty(\widehat{G})$.

 \item For two functions $X,Y:\widehat{G}\rightarrow \mathbb{C}^M$ the notation $X=Y$ means that $X(\xi)=Y(\xi)$ a.e. $\xi\in \widehat{G}$. \end{enumerate}

%\medskip
%%%%%%%%%%%%%%%%%%%%%%%%%%%%%%%%%%%%
%%%%%%%%%%%%%%%%%%%%%%%%%%%%%%%%%%%%%
\section{Discrete convolution operators} \label{LTI1}\label{s1}
%%%%%%%%%%%%%%%%%%%%%%%%%%%%%%%%%%%%

We said that  $\mathcal{A}:\ell^2_{_{N}}(G)\to \ell^2_{_{M}}(G)$  is  a LTI operator if it commutes  with translations, i.e.  $\mathcal{A}\,T_{g}= T_{g} \,\mathcal{A}$, for all $g\in G$. A bounded LTI operator $\mathcal{A}$ can be expressed in the form
\[
\begin{split}
&\big[\mathcal{A} (\mathbf{x})\big]_{1}= \, a_{_{1,1}}\ast x_{_{1}}\, \,\,+\ldots+a_{_{1,N}}\ast x_{_{N}}\\
&\quad \,\,\vdots \qquad \qquad \quad\vdots \qquad\qquad\qquad \, \, \quad\vdots\\
&\big[\mathcal{A} (\mathbf{x})\big]_{_{M}}= a_{_{M,1}}\ast x_{_{1}}+\ldots+a_{_{M,N}}\ast x_{_{N}}
\end{split}
\]
 or in matrix notation
\begin{equation*}
\mathcal{A}(\mathbf{x})= A  \ast \mathbf{x},\quad \mathbf{x}\in \ell^2_{_N}(G),\qquad A=\big[a_{m,n}\big]\in \mathcal{M}_{_{M\times N}}(\ell^2(G)).
\end{equation*}
Indeed, if $\mathcal{A}$ is LTI and bounded,
\begin{equation}\label{conmuta}
[\mathcal{A}(\mathbf{x})]_{m}=\Big[\mathcal{A}\big(\sum_{n=1}^N\sum_{g\in G}x_{n}(g)T_{g}(\delta \mathbf{e}_{n})\big)\Big]_{m}=
\sum_{n=1}^N\sum_{g\in G}x_{n}(g)T_{g}[\mathcal{A}(\delta \mathbf{e}_{n})]_{m}=\sum_{n=1}^N x_n \ast a_{m,n},
\end{equation}
for all $\mathbf{x}\in \ell^2_{_N}(G)$, where $a_{m,n}=\big[\mathcal{A} (\delta \mathbf{e}_{n})\big]_{m}$, $\delta$ is the dirac delta and $\mathbf{e}_{n}$ the  nth column of the identity matrix $\mathbb{I}_{N}$. Reciprocally, a convolution operator $\mathcal{A}:\ell^2_{_{N}}(G)\to \ell^2_{_{M}}(G)$, given by $\mathcal{A}(\mathbf{x})=A\ast \mathbf{x}$, with $A\in  \mathcal{M}_{_{M\times N}}(\ell^2(G))$ commutes with translations.
 In the following theorem we characterize when  it is bounded in terms of
  the Fourier transform of the matrix $A$, which we call, as usual, the transfer matrix of $\mathcal{A}$,
 \[
 \widehat{A}=\big[\, \widehat{a}_{m,n}\, \big] \in \mathcal{M}_{_{M\times N}}(L^2(\widehat{G})). \]

 \medskip
 %\medskip

 To prove the theorem we need
  the Lemma 2.2 of ref. \cite{garcia:19}:

%%%%%%% lema conv
 \begin{lema} \label{conv}
If $x,y\in \ell^2(G)$ and $\widehat{x}\cdot \widehat{y}\in L^2(\widehat{G})$ then
$x\ast y\in \ell^2(G)$   and
 $
 \widehat{x\ast y}=\widehat{x}\cdot \widehat{y}.
 $ \end{lema}

\medskip

 %%%%%%Teorema bessel

   \begin{teo}\label{bessel} Let $A\in  \mathcal{M}_{_{M\times N}}\big(\ell^2(G)\big)$. Then,
  $
  \mathcal{A}: \mathbf{x} \mapsto A  \ast \mathbf{x}
  $
  is a well defined bounded operator from $\ell^2_{_N}(G)$ into $\ell^2_{_M}(G)$ if and only if   $\widehat{A}\in \mathcal{M}_{_{M\times N}}\big(L^\infty(\widehat{G})\big)$.
\end{teo}

\medskip
\noindent{\bf Proof}.   Assume first that $\mathcal{A}$ is a well defined bounded operator  from $\ell^2_{_N}(G)$  into $\ell^2_{_M}(G)$. Denoting by $a_{m,n}$ the $m,n$ entry  of $A$ and $\mathbf{e}_{n}$ the  nth column of  $\mathbb{I}_{N}$, we have that, for any $x\in \ell^2(G)$,
\begin{equation}\label{b1}
\|a_{m,n}\ast x\|_{\ell^2(G)}=\big\|[A\ast x\, \mathbf{e}_{n}]_{m}\big\|_{\ell^2(G)}\le \big\|\mathcal{A}(x\,\mathbf{e}_{n})\big\|_{\ell_M^2(G)}\le \|\mathcal{A}\|  \, \|x\,\mathbf{e}_{n}\|_{\ell_N^2(G)} =
\|\mathcal{A}\|  \, \|x\|_{\ell^2(G)}.
\end{equation}
On the other hand, let us assume that $C$ is  a number such that
$|\widehat{a}_{m,n}(\xi)|>C$ for all $\xi$ in a set $\Omega$ of positive measure.  Since $\widehat{G}$ is compact, $\chi_{_\Omega}$, the characteristic function of $\Omega$,   belongs to $L^2(\widehat{G})$ and then
$x=\mathcal{F}^{-1}\chi_{_\Omega}\in \ell^2(G)$.
By applying Lemma \ref{conv}, we obtain, that for such $C$ and $x$ it is satisfied
\[
\| a_{m,n}\ast x  \|_{\ell^2(G)}= \| \widehat{a_{m,n}\ast x} \|_{L^2(\widehat{G})}=
 \| \widehat{a}_{m,n}\cdot \chi_{_\Omega} \|_{L^2(\widehat{G})}\ge C  \, \| \chi_{_\Omega} \|_{L^2(\widehat{G})}= C \, \| x\|_{\ell^2(G)},
\]
and then, from \eqref{b1}, $C\le \|\mathcal{A}\|$. Hence,  we deduce $|\widehat{a}_{m,n}(\xi)|\le \|\mathcal{A}\|$ a.e. and then $\widehat{a}_{m,n}\in L^\infty(\widehat{G})$.

Reciprocally,  we assume that $\widehat{A}\in \mathcal{M}_{_{M\times N}}\big(L^\infty(\widehat{G})\big)$.  By applying Lemma \ref{conv}, we obtain that $\widehat{a_{m,n}\ast x}= \widehat{a}_{m,n}\cdot \widehat{x}$, for all $x\in \ell^2(G)$. Then, for any $x\in \ell^2(G)$,
\[\begin{split}\|a_{m,n}\ast x\|_{\ell^2(G)}&= \| \widehat{a_{m,n}\ast x} \|_{L^2(\widehat{G})}=
 \| \widehat{a}_{m,n}\cdot \widehat{x} \|_{L^2(\widehat{G})}\le
  \| \widehat{a}_{m,n}\|_{L^\infty(G)} \|x\|_{\ell^2(G)}.
  \end{split}\]
Using this inequality  and \eqref{conmuta}, it follows easily that $\mathcal{A}$ is bounded.  \hfill$\square$

%%%%%%%% Teorema bessel 2
\medskip

\medskip
%%%%%%%%%%%%%%% Remark l^1
\begin{remark}  If $A\in  \mathcal{M}_{_{M\times N}}\big(\ell^1(G)\big)$
then $\widehat{A}\in  \mathcal{M}_{_{M\times N}}\big(C(\widehat{G})\big)\subseteq  \mathcal{M}_{_{M\times N}}\big(L^\infty(\widehat{G})\big)$ and then $\mathcal{A}$ is bounded. The same is true for  $A\in  \mathcal{M}_{_{M\times N}}\big(C^*(G)\big)$, where $C^*(G)=\mathcal{F}^{-1}(C(\widehat{G}))$, the group C*-algebra of $G$ \cite[Section 7.1]{folland:95}.
\end{remark}
\medskip

  \begin{teo} \label{bessel2} Let $\mathcal{A}:\ell^2_{_N}(G)\to \ell^2_{_M}(G)$ a linear operator. Then the following staments are equivalent:
 \begin{itemize}
  \item[(a)] $\mathcal{A}$  commutes with translations and is bounded.
  \item[(b)] There exists $A\in \mathcal{M}_{_{M\times N}}(\ell^2(G))$, such that $\widehat{A} \in \mathcal{M}_{_{M\times N}}\big(L^\infty(\widehat{G})\big)$ and $\mathcal{A}(\mathbf{x})=A\ast \mathbf{x}$ for all $\mathbf{x}\in \ell^2_{_N}(G)$.
  \item[(c)]  There exists $\Lambda \in \mathcal{M}_{_{M\times N}}\big(L^\infty(\widehat{G})\big)$ such that $\widehat{\mathcal{A}(\mathbf{x})}=\Lambda \cdot \widehat{\mathbf{x}}$ for all $\mathbf{x}\in \ell^2_{_N}(G)$.
  \end{itemize}

  In this case, the matrices $A$ and $\Lambda$ in (b) and (c) are unique and $\Lambda=\widehat{A}$.
\end{teo}

\medskip
\noindent{\bf Proof.}

\noindent{\bf (a) $\Leftrightarrow$ (b)} If (a) holds,  from \eqref{conmuta},  $\mathcal{A}(\mathbf{x})=A\ast \mathbf{x}$, $\mathbf{x}\in \ell^2_{_{N}}(G)$, where $A:=[a_{m,n}]$
with $a_{m,n}=\big[\mathcal{A} (\delta \mathbf{e}_{n})\big]_{m}$. From Theorem \ref{bessel}, $\widehat{A}\in \mathcal{M}_{_{M\times N}}\big(L^\infty(\widehat{G})\big)$.
Reciprocally, if (b) holds, $\mathcal{A}$ is bounded from Theorem \ref{bessel} and it commutes with translations since convolution operators do.
\medskip

\noindent{\bf (b) $\Leftrightarrow$ (c)} If (b) holds,   by applying Lemma \ref{conv}, we obtain that  $\widehat{\mathcal{A}(\mathbf{x})}=\widehat{A} \cdot \widehat{\mathbf{x}}$, for $\mathbf{x}\in \ell^2_{_N}(G)$, which proves (c). Reciprocally, we assume that  (c) holds. Since $\widehat{G}$ is compact, $ L^\infty(\widehat{G}) \subseteq L^2(\widehat{G})$. Thus,  $\Lambda\in \mathcal{M}_{_{M\times N}}\big(L^2(\widehat{G})\big)$
and then $A=\mathcal{F}^{-1}(\Lambda)\in\mathcal{M}_{_{M\times N}}\big(\ell^2(G)\big)$. Since the entries $\widehat{a}_{m,n}$ of $\widehat{A}=\Lambda$,  belong to $L^\infty(\widehat{G})$ by using Lemma \ref{conv} we obtain that, for any $\mathbf{x} \in \ell^2_{_N}(G)$,
\[
\widehat{\mathcal{A}(\mathbf{x})}= \Lambda\cdot  \widehat{\mathbf{x}}= \widehat{A}\cdot  \widehat{\mathbf{x}}=
\Big[\sum_{n=1}^N \widehat{a}_{m,n} \cdot \widehat{x}_{n}\Big]_{m=1,\ldots,M}=
\Big[\sum_{n=1}^N \widehat{a_{m,n} \ast x_{n}}\Big]_{m=1,\ldots,M}
\]
and then $\mathcal{A}(\mathbf{x})=\Big[\sum_{n=1}^N a_{m,n} \ast x_{n}\Big]_{m=1,\ldots,M}=A\ast \mathbf{x}$, which proves (b).

\medskip
The matrix $A$ satisfying (b) is unique since,  from \eqref{conmuta}, the entry $m,n$ of such matrix is necessarily given by $a_{m,n}=\big[A\ast  (\delta \mathbf{e}_{n})\big]_{m}$. If $\Lambda$ satisfies
(c) we have proved that
$\mathcal{A}(\mathbf{x})=\mathcal{F}^{-1}(\Lambda)\ast \mathbf{x}$ for all $\mathbf{x}\in \ell_{N}^2(G)$. Therefore, the matrix $\Lambda$ is necessarily $\mathcal{F}(A)$, where $A$ is the unique matrix satisfying (b). \hfill$\square$

\medskip

 For the simplest scalar case $N=M=1$, Theorem \ref{bessel2} can be obtained as a corollary of \cite[Theorems 4.1.1 and 4.3.1]{larsen:71}.

%\medskip

%%%%%%%%%% SUBSECTION Algebra of multipliers
%%%%%%%%%%%%%%%%%%%%%%%%%%%%

\subsection{The algebra $\mathcal{B}_{_{LTI}}\big(\ell^2_{_N}(G)\big)$}

Let us denote by $\mathcal{B}_{_{LTI}}\big(\ell^2_{_N}(G)\big)$ the set of bounded LTI endomorphism of $\ell_N^2(G)$,
\[
\mathcal{B}_{_{LTI}}\big(\ell^2_{_N}(G)\big):=\Big\{\mathcal{A}\in \mathcal{B}\big(\ell^2_{_N}(G)\big): \, \mathcal{A}\, T_{g}=T_{g}\, \mathcal{A}\, \text{ for all }\, g\in G\Big\}.
\]

\medskip

The following lemma proves that $\mathcal{B}_{_{LTI}}\big(\ell^2_{_N}(G)\big)$ is an algebra isomorphic to $\mathcal{M}_{_N}\big(L^\infty(\widehat{G})\big)$. Recall that $\mathcal{M}_{_N}\big(L^\infty(\widehat{G})\big)$ is an involution algebra. The product is the pointwise multiplication and the  involution of a matrix $\Lambda$  is the adjoint matrix $\Lambda^*$.

%%%%%%%%% Teorema algebra
\begin{lema}\label{algebra} $\mathcal{B}_{_{LTI}}\big(\ell^2_{_N}(G)\big)$ is a $*$-subalgebra of $\mathcal{B}\big(\ell^2_{_N}(G)\big)$ and it is $*$-isomorphic to the algebra $\mathcal{M}_{_N}\big(L^\infty(\widehat{G})\big)$. Namely,
\[
\mathcal{B}_{_{LTI}}\big(\ell^2_{_N}(G)\big) \ni  \mathcal{A} \, \longmapsto \,\widehat{A}\in \mathcal{M}_{_N}(L^\infty(\widehat{G}))
\]
is a $*$-isomorphism.
\end{lema}

\noindent{\bf Proof.}  The transform $\mathcal{A} \, \to \,\widehat{A}$ is obviously linear and Theorem \ref{bessel2} proves that it is bijective from $\mathcal{B}_{_{LTI}}\big(\ell^2_{_N}(G)\big)$ onto $\mathcal{M}_{_N}(L^\infty(\widehat{G}))$.
Let $\mathcal{A}, \mathcal{B}\in \mathcal{B}_{_{LTI}}\big(\ell^2_{_N}(G)\big)$.   By applying
Theorem \ref{bessel2} twice  we obtain that
\[\widehat{\mathcal{B}\mathcal{A}(\mathbf{x})}=\widehat{B}\cdot
\widehat{\mathcal{A}(\mathbf{x})}=\widehat{B}\cdot
\widehat{A}\cdot\widehat{\mathbf{x}}\quad \mathbf{x}\in \ell_{N}^2(G).\]
where $\widehat{A},\widehat{B}\in  \mathcal{M}_{_{N}}(L^\infty(\widehat{G}))$. Since $\widehat{B}\cdot \widehat{A}\in  \mathcal{M}_{_{N}}(L^\infty(\widehat{G}))$ and Theorem \ref{bessel2},
$\mathcal{B} \mathcal{A}\in \mathcal{B}_{_{LTI}}\big(\ell^2_{_N}(G)\big)$ and $\mathcal{B} \mathcal{A} \to \widehat{B}\cdot \widehat{A}$. Using Teorema \ref{bessel2} we obtain, that for all $\mathbf{x},\mathbf{y}\in \ell^2_{N}(G)$,
\[
\begin{split}
\langle\widehat{\mathbf{x}},\widehat{\mathcal{A}^*\mathbf{y}}\rangle_{L^2_{N}(\widehat{G})}&=
\langle \mathbf{x},\mathcal{A}^*\mathbf{y}\rangle_{\ell_{N}^2(G)}=
\langle \mathcal{A}\mathbf{x},\mathbf{y}\rangle_{\ell_{N}^2(G)}=
\langle\widehat{\mathcal{A}\mathbf{x}},\widehat{\mathbf{y}}\rangle_{L^2_{N}(\widehat{G})}=
\langle\widehat{A}\cdot \widehat{\mathbf{x}},\widehat{\mathbf{y}}\rangle_{L^2_{N}(\widehat{G})}
\\&=
\langle\widehat{\mathbf{x}},(\widehat{A})^*\cdot \widehat{\mathbf{y}}\rangle_{L^2_{N}(\widehat{G})}
\end{split}
\]
Hence $\widehat{\mathcal{A}^*\mathbf{y}}=(\widehat{A}\,)^*\cdot \widehat{\mathbf{y}}$, for all
$\mathbf{y}\in \ell^2_{N}(G)$. Hence, since $(\widehat{A}\,)^*\in \mathcal{M}_{_N}(L^\infty(\widehat{G}))
$ and Theorem \ref{bessel2} we obtain that
$\mathcal{A}^*\in \mathcal{B}_{_{LTI}}\big(\ell^2_{_N}(G)\big)$ and $\mathcal{A}^* \to (\widehat{A}\,)^*$. \hfill$\square$

\medskip
\begin{remark}\label{r2} Let denote by
\[P(G):=\mathcal{F}^{-1}\big[L^\infty(\widehat{G})\big]=\big\{ \, x \in
\ell^2(G)\, : \, \widehat{x}\in L^\infty(\widehat{G}) \,\big\},\]
to the space of pseudomeasures on $G$, see  \cite[Section 4.2]{larsen:71} and \cite[Section 3.1.8]{benedetto:96}. From Lemma \ref{algebra}, it could be easily proved that the set of matrices $\mathcal{M}_{N}(P(G))$ is an involution algebra where the product is the convolution and the involution $A^*$ of $A$, is not the adjoint of $A$, but
\begin{equation}\label{involution}
A^*=\big[a^*_{m,n}\big]^\top\in  \mathcal{M}_{_{N}}(P(G)),\quad\text{where}\quad a^*_{m,n}(g):=\overline{a_{m,n}(-g)}, \quad g\in G.
\end{equation}
Besides
$
\mathcal{B}_{_{LTI}}\big(\ell^2_{_N}(G)\big) \ni  \mathcal{A} \, \longmapsto \,A \in \mathcal{M}_{_{N}}(P(G))
$
is a $*$-isomorphism.
\end{remark}

%%%%%%%%%%%% The norm
%%%%%%%%%%%%%%%%%%%%%%%%
\subsection{The norm in $\mathcal{B}_{_{LTI}}\big(\ell^2_{_N}(G)\big)$}

The algebra $\mathcal{M}_{_N}\big(L^\infty(\widehat{G})\big)$ is  a C*-algebra, namely it is the matrix algebra on the C*-algebra $L^\infty(\widehat{G})$ \cite[Section II.6.6]{blackadar:06}. The faithful representation of  $L^\infty(\widehat{G})$,  $\sigma:L^\infty(\widehat{G})\rightarrow \mathcal{B}(L^2(\widehat{G}))$, defined by $\sigma(y)x=
y \cdot x$ gives the  faithful representation of $\mathcal{M}_{_N}\big(L^\infty(\widehat{G})\big)$,
\[\gamma:\mathcal{M}_{_N}\big(L^\infty(\widehat{G})\big)\rightarrow \mathcal{B}(L_{_N}^2(\widehat{G})),\quad \gamma(\Lambda)X=
\Lambda \cdot X,\]
which provides an expression for the C*-norm of $\mathcal{M}_{_N}\big(L^\infty(\widehat{G})\big)$,
\begin{equation}\label{norma1}
\|\Lambda\|_{\mathcal{M}_{_N}(L^\infty(\widehat{G}))}= \|\gamma(\Lambda)
\|_{\mathcal{B}(L_{_N}^2(\widehat{G}))}=
\sup_{\|X\|_{L^2_{_{N}}(\widehat{G})}=1} \|\Lambda \cdot X\|_{L^2_{_{N}}(\widehat{G})}\,.
\end{equation}
The algebra $\mathcal{M}_{_N}\big(L^\infty(\widehat{G})\big)$ can also be seen as the algebra $L^\infty(\widehat{G},\mathcal{M}_{_N}(\mathbb{C}))$, that is, the algebra of the essentially bounded functions from $\widehat{G}$ into the C*-algebra $\mathcal{M}_{_N}(\mathbb{C})$. This provides a simpler  expression for the C*-norm of $\mathcal{M}_{_N}\big(L^\infty(\widehat{G})\big)$,
 \begin{equation}\label{norma2}
 \|\Lambda\|_{\mathcal{M}_{_N}(L^\infty(\widehat{G}))}=\esup_{\xi\in \widehat{G}} \| \Lambda(\xi)\|_{2}.
 \end{equation}

\begin{teo}\label{algebra2} $\mathcal{B}_{_{LTI}}\big(\ell^2_{_N}(G)\big)$ is a C*-subalgebra of $\mathcal{B}\big(\ell^2_{_N}(G)\big)$ which is $*$-isometric to the C*-algebra $\mathcal{M}_{_N}\big(L^\infty(\widehat{G})\big)$. Namely
\[
\mathcal{B}_{_{LTI}}\big(\ell^2_{_N}(G)\big) \ni  \mathcal{A} \, \longmapsto \,\widehat{A}\in \mathcal{M}_{_N}(L^\infty(\widehat{G}))
\]
is an isometry $*$-isomorphism (a  C*-isomorphism).
\end{teo}

\noindent{\bf Proof.} By using Theorem \ref{bessel2}, that the Fourier transform is a isometric isomorphism from $\ell^2_{_N}(G)$ onto $L^2_{_N}(\widehat{G})$,  and \eqref{norma1}, we obtain that, for any $  \mathcal{A}\in \mathcal{B}_{_{LTI}}\big(\ell^2_{_N}(G)\big)$,
\[
\begin{split}
\|\mathcal{A}\|&=
\sup_{\|\mathbf{x}\|_{\ell_{_{N}}^2(G)}=1}\|\mathcal{A}(\mathbf{x})\|_{\ell_{_{N}}^2(G)}=
\sup_{\|\mathbf{x}\|_{\ell_{_{N}}^2(G)}=1}\|\widehat{\mathcal{A}(\mathbf{x})}\|_{L^2_{_N}(\widehat{G})}=
\sup_{\|\mathbf{x}\|_{\ell_{_{N}}^2(G)}=1}\|\widehat{A}\cdot \widehat{\mathbf{x}}\|_{L^2_{_N}(\widehat{G})}\\&=
\sup_{\|X\|_{L^2_{_N}(\widehat{G})}=1}\|\widehat{A}\cdot X\|_{L_{_{N}}^2(\widehat{G})}=\|\widehat{A}\,\|_{\mathcal{M}_{_N}(L^\infty(\widehat{G}))}.
\end{split}
\]
Now the theorem follows from Lemma \ref{algebra}.\hfill$\square$

\medskip
 From Theorem \ref{algebra2}, having in mind the expression for the norm \eqref{norma2}, we obtain the following consequence.

\begin{cor}\label{norma}  For any $\mathcal{A}\in \mathcal{B}_{_{LTI}}\big(\ell^2_{_N}(G)\big) $ with transfer matrix $\widehat{A}$, we have
\[ \|\mathcal{A}\|=
\esup_{\xi\in \widehat{G}} \| \widehat{A}(\xi)\|_{2}.\]
\end{cor}

For the case $G=\mathbb{Z}$ this expression for the norm was  given in \cite[Theorem 2.2]{aldroubi:96}.

\medskip
%%%%%%%% Remarks 3

\begin{remark} The norm given by this corollary is difficult to compute. Reference \cite[Theorem 3]{benzi:14} provides the estimacion
$\|\mathcal{A}\|=\|\widehat{A}\,\|_{\mathcal{M}_{N}(L^\infty(\widehat{G}))}\le \|\mathbf{A}\|_{2}$ where $\mathbf{A}=\big[\|\widehat{a}_{m,n}\|_{L^\infty(\widehat{G})}\big]\in \mathcal{M}_{_N}(\mathbb{R})$.
\end{remark}

\begin{remark} The involution algebra $\mathcal{M}_{_N}(P(G))$ defined in Remark \ref{r2} with the norm
$\|A\|=\esup_{\xi\in \widehat{G}} \| \widehat{A}(\xi)\|_{2}
$ is a C*-algebra which is C*-isometric to $\mathcal{B}_{_{LTI}}\big(\ell^2_{_N}(G)\big)$ and to $\mathcal{M}_{_N}\big(L^\infty(\widehat{G})\big)$.
\end{remark}

\subsection{The  invertible elements in $\mathcal{B}_{_{LTI}}\big(\ell^2_{_N}(G)\big)$}

The following theorem characterizes the units of the algebra $\mathcal{B}_{_{LTI}}\big(\ell^2_{_N}(G)\big)$ and provides an expression for the norm of the inverse.

\begin{teo}\label{inversa} Let $\mathcal{A} \in \mathcal{B}_{_{LTI}}\big(\ell^2_{_N}(G)\big)$ with transfer matrix $\widehat{A}$. Then
$\mathcal{A}$ is  invertible in the C*-algebra $\mathcal{B}_{_{LTI}}\big(\ell^2_{_N}(G)\big)$ if and only if $\einf_{\xi\in \widehat{G}} |\det \widehat{A}(\xi)|>0$. In this case
\[
\| \mathcal{A}^{-1}\|= \big(\einf_{\xi\in \widehat{G}} \lambda_{\min} [\widehat{A}(\xi)^*\widehat{A}(\xi)]\big)^{-1/2}.
\]
\end{teo}

\noindent{\bf Proof.} We will prove the equivalent assertion (see Theorem \ref{algebra2}):
$\Lambda\in \mathcal{M}_{_N}(L^\infty(\widehat{G}))$ is  invertible in the C*-algebra $\mathcal{M}_{_N}(L^\infty(\widehat{G}))$
if and only if $\einf_{\xi\in \widehat{G}} |\det \Lambda(\xi)|>0$, and in this case
\[
\| \Lambda^{-1}\|_{\mathcal{M}_{_N}(L^\infty(\widehat{G}))}= \big(\einf_{\xi\in \widehat{G}} \lambda_{\min} [\Lambda(\xi)^*\Lambda(\xi)]\big)^{-1/2}.
\]
If $\einf_{\xi\in \widehat{G}} |\det\Lambda(\xi)|>0$ then there exists the inverse matrix $[\Lambda(\xi)]^{-1}$ a.e. $\xi\in \widehat{G}$. Besides, since $\Lambda\in \mathcal{M}_{_N}(L^\infty(\widehat{G}))$ and $\einf_{\xi\in \widehat{G}} |\det\Lambda(\xi)|>0$, we deduce that    $\Lambda^{-1}\in \mathcal{M}_{_N}(L^\infty(\widehat{G}))$.
    Reciprocally, if  $\Lambda^{-1}\in \mathcal{M}_{_N}(L^\infty(\widehat{G}))$ then $\det \Lambda^{-1} \in
L^\infty(\widehat{G})$, and then, having in mind that $\det \Lambda\det \Lambda^{-1}=1$, we deduce that $\einf_{\xi\in \widehat{G}} |\det \Lambda(\xi)|>0$.

In order to obtain the expression for the norm, note that since
$|\det \Lambda(\xi)|>0$ a.e. we have $\lambda_{\text{min}} [\Lambda(\xi)^*\Lambda(\xi)]>0$ a.e.  Having in mind this fact and \eqref{norma2}, we obtain that
\[
 \begin{split}
 \|\Lambda^{-1}\|_{\mathcal{M}_{_N}(L^\infty(\widehat{G}))}^{-2}&=\big[ \esup_{\xi\in \widehat{G}}\|\Lambda^{-1}(\xi)\|_2\big]^{-2}=  \big[\esup_{\xi\in \widehat{G}} \lambda^{-1}_{\text{min}} [\Lambda(\xi)^*\Lambda(\xi)]\big]^{-1}\\&= \einf_{\xi\in \widehat{G}} \lambda_{\text{min}} [\Lambda(\xi)^*\Lambda(\xi)].\hspace{6cm} \square
 \end{split}
 \]
 Whenever  $\widehat{A}\in  \mathcal{M}_{_{M\times N}}\big(C(\widehat{G})\big)$, the characterization in Theorem \ref{inversa} can be obtained as a corollary of  \cite[Theorem 4]{benzi:14} and Theorem \ref{algebra2}.
%%%%%%%%%%%%%%%%%%%%%%%%%%%%%%%%
%%%%%%%%%%%%%%%%%%%%%%%%%%%%%%%%
\section{Riesz systems generated by actions of abelian groups}\label{s2}

Let $\mathcal{H}$ be a separable complex Hilbert space and $\pi$ a unitary representation of the group $G$ on $\mathcal{H}$, i.e. a homomorphism $\pi:G\rightarrow U(\mathcal{H})$, where  $U(\mathcal{H})$ denotes the group of unitary operators of $\mathcal{H}$, that satisfies
$\pi_{g+g'}=\pi_g\pi_{g'}$, $\pi_{-g}=\pi_g^{-1}=\pi_g^*$, for all $g,g'\in G$.

For a set $\Phi=\{\varphi_{1},\ldots,\varphi_{_N}\}$ of $N$ elements
 of $\mathcal{H}$, we consider the system
\[
\big\{\pi_g\varphi_{n}\big\}_{n\in \mathcal{N},g\in G}\, ,
\]
where $\mathcal{N}=\{1,2,\ldots,N\}$. It is said that $\big\{\pi_g\varphi_{n}\big\}_{n\in \mathcal{N},g\in G}$ is a Riesz sequence of $\mathcal{H}$ when there exist constants $0<\alpha\le\beta<\infty$ such that
\begin{equation}\label{riesz}
\alpha \, \|\mathbf{x}\|^2_{\ell^2_N(G)}\le \Big\| \sum_{n\in \mathcal{N},g\in G} x_{n}(g)\pi_g\varphi_{n} \Big\|_\mathcal{H}^2\le \beta\, \|\mathbf{x}\|^2_{\ell^2_N(G)},
\end{equation}
for all $\mathbf{x}=(x_1,\ldots,x_N)\in \ell^2_N(G)$.
In this case, the system $\big\{\pi_g\varphi_{n}\big\}_{n\in \mathcal{N},g\in G}$ is a Riesz basis for the space \cite[Section 3.6]{ole:03}
 \[
 V_\Phi=\overline{\text{span}} \{\pi_g\varphi_{n}\}_{n\in \mathcal{N},g\in G}=\Big\{\sum_{n\in \mathcal{N},g\in G} x_{n}(g)\pi_{g}\varphi_{n}\, :\,
 (x_{1},\ldots,x_{N})\in \ell_{N}^2(G)\Big\}.
 \]
Note that $V_\Phi$ is invariant by actions   of the group $G$.
When $\mathcal{H}=L^2(\mathbb{R}^d)$, $G=\mathbb{Z}^d$ and $\pi_gf=T_g f=f(\cdot-g)$, these spaces, called shift-invariant spaces, have been widely studied given its importance in wavelets and approximation theory.

\medskip
 The largest $\alpha$ and the smallest $\beta$ satisfying \eqref{riesz} are called the optimal Riesz bounds.
When the right inequality in \eqref{riesz} holds, it is said that
$\big\{\pi_g\varphi_{n}\big\}_{n\in \mathcal{N},g\in G}$ is a Bessel sequence of $\mathcal{H}$ with Bessel bound $\beta$.

\medskip

Let
\begin{equation}\label{A2}
a_{m,n}(g):=\big\langle \, \varphi_{n}\, ,\, \pi_g\varphi_{m}\,\big\rangle_{\mathcal{H}},\qquad A(g):=\big[a_{m,n}(g)\big]\, \in \,\mathcal{M}_{N}(\mathbb{C}),\qquad g\in G.
\end{equation}

\medskip
\noindent For any finite sequences $\mathbf{x}=(x_{1},\ldots,x_{N})$ and $ \mathbf{y}=(y_{1},\ldots,y_{N})$, we have
\begin{equation}
\begin{split}
&\Big\langle \sum_{n\in \mathcal{N},g'\in G} x_{n}(g')\pi_{g'}\varphi_{n},
\sum_{m\in \mathcal{N},g\in G} y_{m}(g)\pi_g\varphi_{m}\Big\rangle_{\mathcal{H}}= \sum_{n,m\in \mathcal{N},g,g'\in G} x_{n}(g')\big\langle \varphi_{n},\pi_{g-g'}\varphi_{m}\big\rangle_{\mathcal{H}}\, \overline{y_{m}(g)}
\\&=  \sum_{n,m\in \mathcal{N},g\in G} (x_{n}\ast a_{m,n})(g)\overline{y_{m}(g)}=
\sum_{m\in \mathcal{N},g\in G} [A\ast \mathbf{x}(g)]_{m}\, \overline{y_{m}(g)}=
\big\langle\,A\ast \mathbf{x}\, ,\, \mathbf{y} \,\big\rangle_{\ell_{N}^2(G)}.\label{AA}
\end{split}
\end{equation}
 Thus, the properties of the system $\big\{\pi_g\varphi_{n}\big\}_{n\in \mathcal{N},g\in G}$  are related to the properties of the convolution operator  $\mathcal{A}(\mathbf{x})=A\ast \mathbf{x}$.

 Note that, whenever equalities \eqref{AA} also hold for any $\mathbf{x},\mathbf{y}\in \ell^2_{_N}(G)$, the operator $\mathcal{A}$ is positive
since in this case  $\big\langle \mathcal{A}(\mathbf{x}), \mathbf{x} \big\rangle_{\ell_{N}^2(G)}=\big\|\sum_{n\in \mathcal{N},q\in G} x_{n}(q)\pi_{q}\varphi_{n}\big\|^2_{\mathcal{H}}\ge 0$ for all $\mathbf{x}\in \ell^2_{_N}(G)$.

In order to use the usual notation in this context, we state the results in terms of the transpose of the transfer matrix (the Gram matrix),
\begin{equation}\label{G}
\mathcal{G}:=\widehat{A}\,^\top,
\end{equation}
which is well defined and belongs to $ \mathcal{M}_{_N}(L^2(\widehat{G}))$ provided that $A$ belongs to $\mathcal{M}_{_N}(\ell^2(G))$.

\medskip
\begin{lema}\label{xAx}
Let  $\mathcal{A}\in \mathcal{B}(\ell_N^2(G))$ be a positive operator.
Then, the operator $\mathcal{A}$ is invertible if and only if $\inf_{\|\mathbf{x}\|=1} \big\langle \mathcal{A}(\mathbf{x}), \mathbf{x} \big\rangle_{\ell_{N}^2(G)}>0$. In this case
$\|\mathcal{A}^{-1}\|=\big(\inf_{\|\mathbf{x}\|=1} \big\langle \mathcal{A}(\mathbf{x}), \mathbf{x} \big\rangle_{\ell_{N}^2(G)}\big)^{-1}.$
 \end{lema}
 \medskip

\noindent{\bf Proof.}  Since $\langle \mathcal{A}(\mathbf{x}), \mathbf{x} \rangle_{\ell_{N}^2(G)}=
\|\mathcal{A}^{1/2}(\mathbf{x})\|_{\ell_{N}^2(G)}^2$,  if $\inf_{\|\mathbf{x}\|=1} \langle \mathcal{A}(\mathbf{x}), \mathbf{x} \rangle_{\ell_{N}^2(G)}>0$, $\mathcal{A}^{1/2}$ is invertible and then $\mathcal{A}$ is invertible. Reciprocally, if $\mathcal{A}$ is invertible then $\mathcal{A}^{1/2}$ is invertible and
\[
\begin{split}
\inf_{\|x\|=1} \langle \mathcal{A}(\mathbf{x}), \mathbf{x} \rangle_{\ell_{N}^2(G)}=\inf_{\|x\|=1}\|\mathcal{A}^{1/2}(\mathbf{x})\|_{\ell_{N}^2(G)}^2
=\|\mathcal{A}^{-1/2}\|^{-2}
=\|\mathcal{A}^{-1/2}\, \mathcal{A}^{-1/2}\|^{-1}=\|\mathcal{A}^{-1}\|^{-1}>0.
\end{split}\]
\hfill$\square$

 \medskip

\begin{teo} \label{R} Let $\varphi_{1}, \varphi_{2},\ldots,\varphi_{N} \in \mathcal{H}$ and $A, \mathcal{G}$ the matrices defined in \eqref{A2} and \eqref{G}. We assume that $A$ belongs to $\mathcal{M}_{_N}(\ell^2(G))$. Then
\begin{itemize}
\item[(a)] The system $\big\{\pi_g\varphi_{n}\big\}_{n\in \mathcal{N},g\in G}$ is a Bessel sequence of $\mathcal{H}$ if and only if $\mathcal{G}\in \mathcal{M}_{_N}(L^\infty(\widehat{G}))$.
In this case, the matrix $\mathcal{G}(\xi)$ is semidefinite positive a.e. $\xi\in \widehat{G}$, and the optimal Bessel bound is $
\esup_{\xi\in \widehat{G}} \lambda_{\max}\,[\mathcal{G}(\xi)].$
\item[(b)]  The system
$\big\{\pi_g\varphi_{n}\big\}_{n\in \mathcal{N},g\in G}$ is a Riesz sequence of
$\mathcal{H}$ if and only if  \[\mathcal{G}\in \mathcal{M}_{_N}(L^\infty(\widehat{G}))\quad \text{and}\quad \einf_{\xi\in \widehat{G}} \det \mathcal{G}(\xi)>0.\] In this case, the optimal Riesz bounds are \[\einf_{\xi\in \widehat{G}} \lambda_{\min}\,[\mathcal{G}(\xi)]\quad \text{and}\quad
\esup_{\xi\in \widehat{G}} \lambda_{\max}\,[\mathcal{G}(\xi)].\]
\end{itemize}
\end{teo}

 \medskip
\noindent{\bf Proof.} Let $\mathcal{A}$ denote the operator defined by $\mathcal{A}(\mathbf{x})= A\ast \mathbf{x}$, $\mathbf{x}\in \ell^2_N(G)$.
\medskip

 \noindent{\bf (a)} We  assume first that $\big\{\pi_g\varphi_{n}\big\}_{n\in \mathcal{N},g\in G}$ is a Bessel sequence of $\mathcal{H}$ with Bessel bound $\beta$. Then, for any $\mathbf{x} \in \ell^2_{_N}(G)$, the series $\sum_{n\in \mathcal{N},q\in G} x_{n}(q)\pi_{q}\varphi_{n}$ converges in $\mathcal{H}$.
Hence, we deduce that the equalities \eqref{AA} hold for any sequences $\mathbf{x},\mathbf{y}\in \ell^2_{_N}(G)$.
For $\mathbf{x}=(x_1,\ldots,x_N)$ let us denote \[f_{\mathbf{x}}=\sum_{n\in \mathcal{N},g\in G} x_{n}(g)\pi_g\varphi_{n}.\] Since $\big\{\pi_g\varphi_{n}\big\}_{n\in \mathcal{N},g\in G}$ is a Bessel sequence with bound $\beta$, we have $\|f_{\mathbf{x}}\|^2_{\mathcal{H}}\le \beta  \|\mathbf{x}\|^2_{\ell_{N}^2(G)}$ for all $\mathbf{x}\in  \ell^2_{N}(G)$.
Then, using  \eqref{AA},
 we obtain that, for any $\mathbf{x},\mathbf{y} \in  \ell_{N}^2(G)$,
\[
\big|\langle \mathcal{A}(\mathbf{x}), \mathbf{y} \rangle_{\ell_{N}^2(G)}\big|=
\big|\langle f_{\mathbf{x}}, f_{\mathbf{y}} \rangle_{\mathcal{H}}\big|\le \|f_{\mathbf{x}}\|_{\mathcal{H}}
 \|f_{\mathbf{y}}\|_{\mathcal{H}} \le \beta \,\|\mathbf{x}\|_{\ell_{N}^2(G)}\, \|\mathbf{y}\|_{\ell_{N}^2(G)}.
\]
Hence, the sesquilinear functional $(\mathbf{x},\mathbf{y})\mapsto\langle \mathcal{A}(\mathbf{x}), \mathbf{y} \rangle_{\ell_{N}^2(G)}$ is bounded. Then, the convolution operator $\mathcal{A}$ belongs to $\mathcal{B}(\ell_{N}^2(G))$. Then, from
Theorem \ref{bessel}, the matrix $\mathcal{G}=\widehat{A}\,^\top$ belongs to  $\mathcal{M}_{_N}(L^\infty(\widehat{G}))$.

Reciprocally, assume now that $\mathcal{G}\in \mathcal{M}_{_N}(L^\infty(\widehat{G}))$. Then, from Theorem \ref{bessel},  the operator $\mathcal{A}$ belongs to $\mathcal{B}(\ell_{N}^2(G))$. By using \eqref{AA}, we obtain that for any finite sequence $\mathbf{x}$,
\[\Big\|\sum_{n\in \mathcal{N},g\in G} x_{n}(g)\pi_g\varphi_{n}\Big\|^2_{\mathcal{H}}=  \big\langle \mathcal{A}(\mathbf{x}), \mathbf{x} \big\rangle_{\ell_{N}^2(G)}\le \, \|\mathcal{A}\|\,  \|\mathbf{x}\|^2_{\ell_{N}^2(G)}.\]
Then, $\big\{\pi_g\varphi_{n}\big\}_{n\in \mathcal{N},g\in G}$ is a Bessel sequence of $\mathcal{H}$ \cite[Theorem 3.6.6]{ole:03}.

\medskip
Assume now that the equivalent conditions of (a) hold. Hence,  \eqref{AA} holds for any sequences $\mathbf{x},\mathbf{y}\in \ell^2_{N}(G)$. Then $\mathcal{A}$ is a positive operator, and thus a positive element of the C*-algebra $B(\ell_{_N}^2(G))$. Hence, it is  a positive element of the C*-subalgebra $B_{LTI}(\ell_{_N}^2(G))$ \cite[II.3.1]{blackadar:06}. Using  Theorem \ref{algebra2}, we obtain that
there exists  $\Lambda\in \mathcal{M}_{_N}(L^\infty(\widehat{G}))$ such that $\widehat{A}=\Lambda^*\Lambda$, and thus the matrix $\widehat{A}(\xi)=\mathcal{G}(\xi)^\top$  is semidefinite positive  a.e. $\xi\in \widehat{G}$. Then, from Corollary \ref{norma} and \eqref{AA},
\[
\begin{split}
\esup_{\xi\in \widehat{G}} \lambda_{\max} \, [\mathcal{G}(\xi)]
&=
\esup_{\xi\in \widehat{G}} \|\widehat{A}(\xi)\|_{2}=
\|\mathcal{A}\|=
 \sup_{\|\mathbf{x}\|_{\ell_{N}^2(G)}=1}\big\langle \mathcal{A}(\mathbf{x}), \mathbf{x} \big\rangle_{\ell_{N}^2(G)}\\&= \sup_{\|\mathbf{x}\|_{\ell_{N}^2(G)}=1} \Big\|\sum_{n\in \mathcal{N},q\in G} x_{n}(q)\pi_g\varphi_{n}\Big\|^2_{_{\mathcal{H}}},
\end{split}
\]
which proves (a).

\medskip
 \noindent{\bf (b)} Since (a), we just have to prove that $ \einf_{\xi\in \widehat{G}} \det \mathcal{G}(\xi)>0$ if and only if the left inequality in \eqref{riesz} holds, and that  the lower optimal Riesz bound is   $ \einf_{\xi\in \widehat{G}} \lambda_{\min} (\mathcal{G}(\xi))$. Besides, we have proved that in any of the hypotheses in (b), the equalities \eqref{AA} hold for any $\mathbf{x},\mathbf{y}\in \ell_{_N}^2(G)$, $\mathcal{A}$ is positive operator, and $\widehat{A}(\xi)$ and $\mathcal{G}(\xi)$ are semidefinite positive a.e. $\xi\in \widehat{G}$.

 Assume first that $ \einf_{\xi\in \widehat{G}} \det \mathcal{G}(\xi)>0$.
 From  Theorem \ref{inversa}, the operator $\mathcal{A}$ is invertible, and \[0<\|\mathcal{A}^{-1}\|^{-2}=\einf_{\xi\in \widehat{G}}\, \lambda_{\min}\, [\widehat{A}(\xi)^*\widehat{A}(\xi)]=  \einf_{\xi\in \widehat{G}} \lambda^2_{\min} [\widehat{A}(\xi)]=\einf_{\xi\in \widehat{G}} \lambda^2_{\min}\, [\mathcal{G}(\xi)].\]
 Hence, using \eqref{AA} and  Lemma \ref{xAx},  we obtain
\[
\inf_{\|x\|=1} \Big\|\sum_{n\in \mathcal{N},q\in G} x_{n}(q)\pi_g\varphi_{n}\Big\|^2_{_{\mathcal{H}}}=
\inf_{\|x\|=1} \big\langle \mathcal{A}(\mathbf{x}), \mathbf{x} \big\rangle_{\ell_{N}^2(G)}=
 \|\mathcal{A}^{-1}\|^{-1}=
\einf_{\xi\in \widehat{G}} \lambda_{\min} \, [\mathcal{G}(\xi)]>0.
\]
Therefore  $\{\pi_g(\varphi_{m})\}_{n\in \mathcal{N},g\in G}$ is a Riesz sequence of $\mathcal{H}$ and the optimal lower Riesz bound is $\einf_{\xi\in \widehat{G}} \lambda_{\min} \mathcal{G}(\xi).$

To prove the reciprocal, assume now that $\big\{\pi_g\varphi_{n}\big\}_{n\in \mathcal{N},g\in G}$ is a Riesz sequence of $\mathcal{H}$. Then, from \eqref{AA},  we have that
\[\inf_{\|x\|=1} \langle \mathcal{A}(\mathbf{x}), \mathbf{x} \rangle_{\ell_{N}^2(G)}=\inf_{\|x\|=1} \Big\|\sum_{n\in \mathcal{N},q\in G} x_{n}(q)\pi_g\varphi_{n}\Big\|^2_{_{\mathcal{H}}}>0.\] Hence, from Lemma \ref{xAx}, we obtain that the operator $\mathcal{A}$ is invertible. Hence, from Theorem \ref{inversa}, $\einf_{\xi\in \widehat{G}} \det \mathcal{G}(\xi)=\einf_{\xi\in \widehat{G}} \det\widehat{A}(\xi)>0$.
\hfill $\square$

\medskip
\medskip

For the classical
shift-invariant systems,
 $\big\{\pi_g\varphi_{n}\big\}_{n\in \mathcal{N},g\in G}=
\big\{T_{g}\varphi_{n}\big\}_{n\in \mathcal{N},g\in \mathbb{Z}^d}\subset L^2(\mathbb{R}^d)$,   the result of this theorem is very well known (see for example \cite{jia:91,ron:95}), and it has many applications in wavelet theory and approximation theory.  It is given  usually in terms of  \[\sum_{\alpha\in \mathbb{Z}^d}
\widehat{\varphi}_{n}(\xi+2\pi \alpha)\overline{\widehat{\varphi}_{m}(\xi+2\pi \alpha)}\] which is equal to $\widehat{a}_{m,n}(e^{i\xi})$
under appropiate conditions, see \cite[eq. 4.1]{hernandez:10} and \cite[Thm. 3.2]{jia:91}.

\medskip

The result given in Theorem \ref{R}  is related to many of the generalizations previously considered: For the systems $\big\{\pi_g\varphi_{n}\big\}_{n\in \mathcal{N},g\in G}=
\big\{U_{1}^{g_{1}}\cdots U_{d}^{g_{d}}\varphi_{n}\big\}_{n\in \mathcal{N},(g_{1},\ldots,g_{d})\in \mathbb{Z}^d}$, where $U_{1},\ldots,U_d$ are unitary operators of $\mathcal{H}$, the result was given in \cite{aldroubi:96,lee:93}; For the systems $\big\{T_{g} \varphi_{n}\big\}_{n\in \mathcal{N},g\in G}$, where the set of generators $\{\varphi_n\}_{n\in \mathcal{N}}\subset L^2(S)$  can be countable, $S$  is a  LCA group, and $G$ is a discrete subgroup of
$S$ such that $S/G$ is compact,  the corresponding result was given in \cite{cabrelli:10}; For the systems $\big\{\pi_{g} \varphi\big\}_{g\in G}$, where the representation $\pi$ satisfies the, so called, dual integrability condition,  the result was given in \cite{hernandez:10}, see \cite{barbieri:15} for the non abelian case with a countable set of generators.

\medskip
\medskip

\medskip
\noindent{\bf Acknowledgments:} The author wishes to thank Antonio Garc\'ia and Miguel Angel Hern\'andez Medina for the stimulating conversations on this work, their suggestions and constructive comments.

%\medskip


\begin{thebibliography}{10}


\bibitem{aldroubi:96}
A. Aldroubi.
\newblock Oblique proyections in atomic spaces.
\newblock {\em Proc. Amer. Math. Soc.}, 124 (1996), pp. 2051--2060.


\bibitem{barbieri:15}
D. Barbieri, E. Hern\'andez and J. Parcet.
\newblock Riesz and frame systems generated by unitary actions of discrete groups.
\newblock {\em Appl. Comput. Harmon. Anal.}, 39 (2015), pp. 369--399.

\bibitem{barbieri:18}
D. Barbieri, E. Hern\'andez and V. Paternostro.
\newblock Invariant spaces under unitary representations of discrete groups.
\newblock { arXiv: 1811.02993}, 2018.


\bibitem{benedetto:96}
J.J.~Benedetto.
\textit{ Harmonic Analysis and Applications}.
CRC Press, Boca Raton FL, 1996.



\bibitem{benzi:14}
M.~Benzi and P.~Boito.
\newblock Decay properties for functions of matrices over C*-algebras.
\newblock {\em Linear Algebra Appl.}, 456 (2014), pp. 174--198.

\bibitem{blackadar:06}
B.~Blackadar.
\newblock {\em Operator Algebras: Theory of C*-Algebras and von Neumann Algebras}.
\newblock Encyclopaedia Math. Sci., Vol. 122, Springer, Berlin, 2006.

\bibitem{cabrelli:10}
C.~Cabrelli and V.~Paternostro.
\newblock Shift-invariant spaces on LCA groups.
 \newblock {\em J. Funct. Anal.}, 258 (2010), pp. 2034--2059.

\bibitem{ole:03}
O.~Christensen.
\newblock {\em An Introduction to Frames and Riesz Bases}.
\newblock Birkh{\"a}user, Boston, 2003.


\bibitem{folland:95}
G.B.~Folland.
\newblock {\em A Course in Abstract Harmonic Analysis}.
\newblock CRC Press, Boca Raton FL, 1995.




\bibitem{garcia:19}
A.G. Garc\'{\i}a and G.~P\'erez-Villal\'on.
\newblock Riesz bases associated with regular representations of semidirect product groups.
\newblock {\em Banach J. Math. Anal.} To appear.

\bibitem{garcia:19b} A.G. Garc\'{\i}a, M.A Hern\'andez-Medina and G.~P\'erez-Villal\'on.
\newblock Convolution systems on discrete abelian groups as
a unifying strategy in sampling theory. {\em Preprint}.


\bibitem{lee:93}
T.N. Goodman, S.L. Lee, and W.S. Tang.
\newblock Wavelet bases for a set of commuting unitary operators. \newblock {\em Adv. Comput. Math.}, 1 (1993), pp. 109-126.


\bibitem{hernandez:10}
E.~Hern\'andez, H.~Sikic, G.~Weiss and E.~Wilson.
\newblock Cyclic subspaces for unitary representations of LCA groups; generalized Zak transform.
\newblock {\em Colloq. Math.}, 118 (2010), pp. 313--332.



\bibitem{jia:91}
R.Q. Jia and C.A. Micchelli.
\newblock Using the refinement equations for the construction of pre-waveles
  {II}: Powers of two.
\newblock {\em In Curves and Surfaces.}
\newblock {\em Academic Press, Boston}, pp. 209--246, 1991.


\bibitem{kailath:80}
T.~Kailath.
\newblock {\em Linear Systems}.
\newblock Prentice  Hall, Berlin, 1980.


\bibitem{larsen:71}
R.~Larsen.
\newblock {\em An Introduction to the Theory of Multipliers}.
\newblock Springer-Verlag, 1971.

\bibitem{ron:95}
A. Ron and Z. Shen.
\newblock Frames and stable bases for shift-invariant subspaces of $L^2(\mathbb{R}^d)$.
\newblock {\em Canad. M. Math.}, 47 (1995), pp. 1051-1094.

\bibitem{rudin4}
W.~Rudin.
\newblock {\em Fourier Analysis on Groups}.
\newblock Interscience Publishers, NewYork - London, 1962.


\end{thebibliography}
\end{document}